\documentclass[12pt]{amsart}
\usepackage{geometry} 
\newtheorem{theorem}{Theorem}[section]
\newtheorem{prop}{Proposition}[section]
\newtheorem{cor}{Corollary}[section]
\newtheorem{rem}{Remark}[section]
\newtheorem{lem}{Lemma}[section]
\newtheorem{defn}{Definition}[section]
\newtheorem{ex}{Example}[section]
\newcommand{\bz}{\mathbb{Z}}
\newcommand{\bq}{\mathbb{Q}}
\newcommand{\bn}{\mathbf{NS}}
\newcommand{\bP}{\mathbf{Pic}}
\newcommand{\wa}{\widehat{A}}
\newcommand{\wb}{\widehat{B}}
\newcommand{\fm}{Fourier-Mukai~}
\newcommand{\bphi}{\mathbf{\Phi}}
\newcommand{\dd}{\mathbf{dim}}
\newcommand{\mo}{\mathcal{O}}
\newcommand{\bd}{\mathbf{det}}
\newcommand{\bext}{\mathbf{Ext}}
\newcommand{\mP}{\mathcal{P}}
\title{Fourier-Mukai Partners of Abelian varieties}
\author{Anningzhe Gao}
\date{} 

\begin{document}

\maketitle
\tableofcontents

\section{Introduction}
In this small note we will consider the Fourier-Mukai partners of an abelian variety over a field of character 0. 
We use $k$ to denote an algebraically closed field of char 0. By an abelian variety over $k$ we mean a proper integral group scheme over $k$ (Since it is over $k$, so it is smooth and projective automatically). Given $a\in A(k)$, we use $T_a$ to denote the translation on $A$ by $a$, so $T_a$ is an automorphism of $A$. Let $X$ be a smooth projective variety over $k$, we use $D^b(X)$ to denote the bounded derived category of coherent sheaves on $X$. Given an abelian variety $A$, we know the Neron-Severi group $\bn(A)=\bP(A)/\bP^0(A)$ is torsion-free. For every element $\mu=\frac{L}{l}\in \bn(A)\otimes_{\bz}\bq$, define 
$$A_\mu=Im(A\to A\times \wa)$$
$$a\to(l_A(a),h_L(a))$$
here $l_A$ is the endomorphism of $A$ defined by multiplied by $l$, and $h_L:A\to\wa$, $a\to T_a^*L\otimes L^{-1}$. Then $A_\mu$ is a sub abelian variety of $A\times \wa$. The following simple result can be concluded from Orlov's result \cite{orlov2002derived} directly. 
\begin{theorem}
Given two abelian varieties $A,B$ over $k$, then $D^b(A)$ is equivalent to $D^b(B)$ as triangulated categories if and only if $\hat{B}\cong A_\mu$ for some $\mu\in\bn(A)\otimes_\bz\bq$ 
\end{theorem}

This theorem will be proved in section 4.

From Orlov's theorem  we know that to give a Fourier-Mukai equivalence between abelian varieties it suffices to consider the kernel for semi-homogeneous vector bundles. That means if we have a derived equivalence between two abelian varieties then we have an equivalence (different with the given one) which is induced by a semi-homogenous vector bundle. In section 5 we will consider the following question: Given a semi-homogenous vector bundle, in which case it will induce an equivalence. More precisely, we have (notations will be given in Section 5)

\begin{theorem}\label{main1}
Let $A,B$ be two abelian varieties. 

(a) If there exists a simple semi-homogeneous vector bundle $E$ on $A\times B$ such that the \fm transform $\bphi_E$ is an equivalence, then:

(1) $\wa\cong B_{\delta(E|_{e_A\times B})}$.

(2) Let $\mu$ be the slope of $E$, then $(A\times B)_\mu$ gives an isomorphism between $A\times\hat{A}$ and $B\times \hat{B}$. This is due to Orlov, but we will give another description by using semi-homogenous vector bundles directly.

(3) Let $A\cong\hat{B_{\delta(E|_{e_A\times B})}}\to\wb$ be the canonical morphism, then $$\bd(E)\cong q_1^*\bd(E|_{A\times e_B})\otimes(\pi\times 1)^*\mP\otimes q_2^*\bd(E|_{e_A\times B})$$

(b) Let $E_B$ be a simple semi-homogeneous vector bundle on $B$ with rank $l$ and slope $\delta$. Let $\pi:B\to B_\delta$ be the natural map. The necessary and sufficient condition for the existence of the equivalence $\bphi_E$ with $\bphi_E(\mo_{e_A})=E_B$ is there exists a line bundle $N$ on $\hat{B_\delta}$ such that $K(N)\cap\hat{B_\delta}_l=Ker(\hat{\pi})$.
\end{theorem}

The organization of this paper is as follows: In section 2 we will discuss some basic facts about Fourier-Mukai transforms, especially in the case of abelian varieties. In section 3 we consider the semi-homogeneous vector bundles on abelian varieties. In section 4 we use the tools from section 2 and section 3 to prove Theorem 1.1. In section 5 we will discuss in the special case when the \fm partner is a vector bundle and discuss Theorem \ref{main1}.
\section{Fourier-Mukai Transforms}
We refer to \cite{orlov2002derived},\cite{orlov2003derived} and \cite{huybrechts2006fourier} for details. For the basic properties of abelian varieties, we refer to \cite{milne1986abelian} and \cite{mumford1974abelian}.

Let $X$ and $Y$ be two smooth projective varieties over $k$, $\dd X=n$, $\dd Y=m$, $k$ is algebraically closed with char 0. Recall that if $X$ is smooth and projective, then any complex $F\in D^b(X)$ is perfect.

Let $P$ be a perfect complex in $D^b(X\times Y)$. Then we can define the \fm transform $$\bphi_P:D^b(X)\to D^b(Y)$$
$$\bphi_P(F)=R\pi_{Y*}(P\otimes \pi_X^*(F))$$
here $\pi_X:X\times Y\to X$ and $\pi_Y:X\times Y\to Y$ are two projections. 

Let $P^\vee=\mathcal{HOM}_{X\times Y}(P,\mathcal{O}_{X\times Y})$. Here $\mathcal{HOM}$ is the sheaf hom. $P^\vee$ can be viewed as a complex on $Y\times X$. Denote $\omega_X$ (resp. $\omega_Y$) is the canonical line bundle on $X$ (resp. $Y$). And let $H=P^\vee\otimes\pi_Y^*\omega_Y[m]$, $G=P^\vee\otimes\pi_X^*\omega_X[n]$, then we have the \fm transform $\bphi_H$ (resp. $\bphi_G$) defined by $H$ (resp. G) from $D^b(Y)$ to $D^b(X)$. Then we have
\begin{prop}(\cite{bridgeland1999equivalences})\label{adjoint}
The functor $\bphi_H$ (resp. $\bphi_G$) is left adjoint (resp. right adjoint) to $\bphi_P$.
\end{prop}

In \cite{orlov2003derived}, Orlov showed any triangulated equivalence between $D^b(X)$ and $D^b(Y)$ is a \fm transform.
\begin{theorem}($\mathbf{Orlov}$)\label{FM Equiv}
Let $X$, $Y$ be two smooth projective varieties over $k$, suppose $F:D^b(X)\to D^b(Y)$ is an equivalence between triangulated categories, then there exists an object $P\in D^b(X\times Y)$, unique up to isomorphisms, such that $F\cong\bphi_P$.
\end{theorem}
In the case of the theorem, we call $P$ the kernel of the \fm transform.

If $\bphi_P$ is an equivalence, then the two adjunction maps 
$$\epsilon:\bphi_H\bphi_P\to id_{D^b(X)}$$
$$\eta:id_{D^b(Y)}\to\bphi_P\bphi_G$$
are isomorphisms, so $\bphi_H\cong\bphi_G$ as the inverse of $\bphi_P$. Then by the uniqueness of the kernel, we have
\begin{prop}(\cite{bridgeland1999equivalences})\label{quasi-inverse}
Let $\bphi_P:D^b(X)\to D^b(Y)$ be a \fm equivalence, then we have $P^\vee\otimes\pi_X^*\omega_X[n]\cong P^\vee\otimes\pi_Y^*\omega_Y[m]$.
\end{prop}

Then we discuss a theorem which can be considered as the "converse" version of Prop \ref{quasi-inverse}. We first establish a lemma which consider the fully faithfulness of a \fm transform.

Let $\bphi_P:D^b(X)\to D^b(Y)$ to be a \fm transform (not necessary an equivalence). For every closed point $x\in X$, let $\mo_x$ be the structure sheaf of the closed point, i.e. the skyscraper sheaf $k(x)$. 
\begin{lem}(Bondal, Orlov)\label{fully-faithful}
The functor $\bphi_P:D^b(X)\to D^b(Y)$ is fully faithful if and only if for any two closed point $x_1,x_2\in X$ one has
$$Hom(\bphi_P(\mo_{x_1}),\bphi_P(\mo_{x_2})[i])=0$$ if $x_1\neq x_2$ or $i<0$ or $i>\dd X$
$$Hom(\bphi_P(\mo_{x_1}),\bphi_P(\mo_{x_2})[i])=k$$ if $x_1=x_2$ and $i=0$
\end{lem}

The following theorem will be used in the proof of Theorem 1.1, which in some sense is the converse version of Prop. \ref{quasi-inverse}
\begin{theorem}(\cite{huybrechts2006fourier})\label{equiv}
Let $\bphi_P:D^b(X)\to D^b(Y)$ be a fully faithful \fm transform between two smooth projective varieties, then $\bphi_P$ is an equivalence if and only if 
$$\dd X=\dd Y$$
$$P\otimes\pi_X^*\omega_X=P\otimes\pi_Y^*\omega_Y$$
\end{theorem}

This theorem is really important in our proof since if we are considering two smooth projective varieties with trivial canonical line bundles (which is the case for abelian varieties!) then the second equation is automatically satisfied. 

Then we put ourselves on the case when $X$ and $Y$ are abelian varieties. We use $A$ and $B$ to stand for abelian varieties. The essential theorem we will use is the following:

\begin{theorem}(Orlov \cite{orlov2002derived})\label{isometry}
If two abelian varieties $A$ and $B$ are derived equivalent, then we have an isomorphism
$$f:A\times\wa\to B\times \wb$$
and in this case $f(a,\alpha)=(b,\beta)$ if and only if
$$\bphi_P(L_\alpha\otimes T_a^*(F))=L_\beta\otimes T_b^*(\bphi_P(F))$$
for any complex $F\in D^b(A)$. Here $L_\alpha$ (resp. $L_\beta$) is the line bundle corresponding to the closed point $\alpha\in\wa$ (resp. $\beta\in\wb$) 
\end{theorem}

With this theorem, we can prove the following interesting corollary.

\begin{cor}
Let $\bphi_P:D^b(A)\to D^b(B)$ be a \fm equivalence. Let $e\in A$ be the closed point corresponding to the unit in $A$. If $\bphi_P(\mo_e)=\mo_b$ for some $b\in B$, then $A$ is isomorphic to $B$. If $\bphi_P(\mo_e)=L$ for some $L$ to be a line bundle, then $A$ is isomorphic to $\wb$.
\end{cor}
\begin{proof}
If $\bphi_P(\mo_e)=L$, then compose $\bphi_P$ with $\otimes L^{-1}:D^b(B)\to D^b(B)$, then the new equivalence, which is also a \fm transform by theorem \ref{FM Equiv}, sends $\mo_e$ to $\mo_B$, then compose with the \fm equivalence $D^b(\wb)\to D^b(B)$ defined by the Poincare bundle on $\wb\times B$, $\mo_e$ will map to $\mo_{\widehat{e}}$, where $\widehat{e}\in\wb$ is the unit of $\wb$. So we just need to prove the first statement.

Let $\bphi_P(\mo_e)=\mo_b$, compose with $T_{-b}^*:D^b(B)\to D^b(B)$, we may assume $b=e'\in B$ to be the unit of $B$. By \cite{favero2012reconstruction}, we just need to prove for any $a\in A$, $\bphi_P(\mo_a)=\mo_b$ for some $b\in B$. Consider the isomorphism $f:A\times\wa\to B\times\wb$ in theorem \ref{isometry}, let $f(a,e_{\wa})=(b,\beta)$, where $e_{\wa}$ is the unit of $\wa$. Then by theorem \ref{isometry}, $$\bphi_P(\mo_a)=\bphi_P(T_a^*(\mo_e))=L_\beta\otimes T_b^*(\bphi_{P}(\mo_e))=L_\beta\otimes T_b^*\mo_{e'}=\mo_b$$
so by \cite{favero2012reconstruction}, $A$ is isomorphic to $B$ (Note that abelian varieties are irreducible non-singular, hence divisorial).
\end{proof}

\section{Semi-homogeneous Vector Bundles}
To give more details of the \fm kernels between abelian varieties, we need to consider semi-homogeneous vector bundles on abelian varieties. Most of material in this section is from \cite{mukai1978semi} and \cite{orlov2002derived}.

\begin{defn}\label{semi-homo}
A vector bundle $M$ on an abelian variety $A$ is called semi-homogeneous if for any $a\in A$, we have $$T^*_a(M)\cong M\otimes L$$ for some $L\in\bP(A)$.
\end{defn}

Let $p_1:A\times \wa\to A$ and $p_2:A\times\wa\to\wa$ be two projections, and let $\mathcal{P}$ be the Poincare bundle. We have 

\begin{defn}
Let $M$ be a simple semi-homogeneous vector bundle on $A$, define a group scheme $\Sigma(M)$ to be the maximal subgroup scheme of $\wa$ such that $p_1^*M\otimes\mathcal{P}$ is isomorphic to $p_1^*M$ over $A\times\Sigma(M)$.
\end{defn}

In char 0 case, $\Sigma(M)$ is reduced, and finite over $k$. It is coincide with $\{\alpha\in\wa,M\otimes L_\alpha\cong M\}$.

Semi-homogeneous vector bundles have the following nice property:
\begin{prop}\label{isogeny}
Let $\pi:A\to B$ be an etale isogeny, and $M$ a semi-homogeneous vector bundle on $A$, then $\pi_*(M)$ is a semi-homogeneous vector bundle on $B$.
\end{prop}
We then will define the $slope$ of a semi-homogeneous vector bundle.

\begin{defn}\label{slope}
Let $M$ be a semi-homogeneous vector bundle. define its slope to be $\frac{\bd(M)}{rank(M)}\in\bn(A)\otimes_\bz\bq$, and we denoted by $\delta(M)$.
\end{defn}

We collect some facts about the semi-homogeneous vectors bundle here. All the details can be found in \cite{mukai1978semi}.

We fix an element $\delta\in\bn(A)\otimes_\bz\bq$, then we have
\begin{prop}\label{im}
(1) There exists a simple semi-homogeneous vector bundle $E$ with slope $\delta$. And if $E$ and $E'$ are two simple semi-homogeneous vector bundle of slope $\delta$, we have $E'\cong E\otimes L$ for some $L\in\bP^0(A)$.

(2) If $E=F\oplus G$ is of slope $\delta$, then $F,G$ have slope $\delta$.

(3) Let $E,E'$ be two simple semi-homogeneous vector bundles of the same slope, if $E$ is not isomorphic to $E'$, then
$$\bext^i(E,E')=\bext^i(E',E)=0$$ for all $i$.
\end{prop}

\begin{defn}
Let $\delta=\frac{L}{l}\in\bn(A)\otimes_\bz\bq$, define $A_\delta=Im(A\to A\times \wa)$ where the morphism is defined by $a\to(la,T_a^*L\otimes L^{-1})$.
\end{defn}

\begin{theorem}(\cite{orlov2002derived})\label{semi}
For a simple semi-homogeneous vector bundle $E$ on $A$ with slope $\mu$, $T_a^*E\cong E\otimes L_\alpha$ if and only if $(a,\alpha)\in A_\mu$.
\end{theorem}

The following is in \cite{orlov2002derived} 
\begin{theorem}(Orlov)\label{Iso-abelian}
If $A$ and $B$ are derived equivalent (that means $D^b(A)$ is equivalent to $D^b(B)$), then there exists $E$ a semi-homogeneous vector bundle on $A\times B$ such that the \fm transform $\bphi_E$  is an equivalence.
\end{theorem}

\section{Some Applications}
In this section we will give an estimation of the numbers of \fm partners of abelian varieties (Here we only consider the \fm partner also an abelian variety). 

For any line bundle $L$ on an abelian variety $A$, $\dim A=g$, let $K(L)$ be the kernel of the morphism $h_L$ defined by $L$, which is $h_L(a)=T_a^*L\otimes L^{-1}$, we know deg$h_L$=(deg$L)^2$ if $L$ is ample. Mumford proved in \cite{mumford1966equations} that if $L$ is non-degenerated (which means $K(L)$ is finite over $k$), then there exists an elementary divisor $\delta=(d_1,d_2,...,d_g)$ (an elementary divisor is an ordered set $(d_1,...,d_g)$ with $d_{i+1}|d_i$), with $K(L)\cong \oplus(\bz/d_i\bz)^2$. This isomorphism is not canonical, actually $K(L)\cong(\oplus\bz/d_i\bz)\oplus(\oplus\mu_{d_i})$ naturally.

We first need a small lemma, here char $k=0$.

\begin{lem}\label{pic}
Let $A$ and $B$ be two abelian varieties over $k$ and $f:A\to B$ is an isogeny, then $rank(\bn(A))=rank(\bn(B))$
\end{lem}
\begin{proof}
We just need to prove the map $f^*:\bn(B)\to\bn(A)$ is injective. Now let's prove $f^*$ is injective. Suppose $f^*(L)=0$ in $\bn(A)$, then $f^*(L)=M\in\bP^0(A)$, since the dual $f^*:\wb\to \wa$ is also an isogeny, hence surjective, then we may choose $L'\in\bP^0(B)$ such that $f^*L'=M$, so we may replace $L$ by $L\otimes L'^{-1}$ and assume $f^*(L)=\mo_B$. Then we have $f_*f^*(L)=f_*(\mo_A)=\oplus_{\beta\in Ker(f^*:\wb\to\wa)}L_\beta$, $L_\beta$ is the line bundle on $B$ corresponding to the point $\beta\in\wb$. By \ref{isogeny}, $f_*f^*L=\oplus_{\beta\in Ker(f^*:\wb\to\wa)}L_\beta$ is semi-homogeneous of slope $0\in\bn(B)\otimes_\bz\bq$, but we also have $f_*f^*L=L\otimes(\oplus_{\beta\in Ker(f^*:\wb\to\wa)}L_\beta)$, so $L$ is a direct summand of $f_*f^*L$, so $L$ is of slope $0$, so $L=0$ in $\bn(B)$, that means $f^*:\bn(B)\to\bn(A)$ is injective. With the discussion in the beginning of the proof, $rank(\bn(A))=rank(\bn(B))$.
\end{proof}

Since any two derived equivalent abelian varieties are isogeneous, then any two derived equivalent abelian varieties have the same Picard number.

Let's prove the essential theorem of this section.

\begin{theorem}\label{essential}
Let $A$ and $B$ be two abelian varieties, let $\bphi_E:D^b(A)\to D^b(B)$ be a \fm equivalence where $E$ is a semi-homogeneous vector bundle. Define $E_B=E|_{e_A\times B}$ ($e_A$ is the unit of $A$), then $\wa\cong B_{\delta(E_B)}$
\end{theorem}
\begin{proof}
Let $\Phi:A\times\wa\to B\times\wb$ to be the isometry defined by $\bphi_E$. Then by theorem \ref{isometry}, we have $\Phi(a,\alpha)=(b,\beta)$ if and only if $$\bphi_E(L_\alpha\otimes T_a^*(F))=L_\beta\otimes T_b^*(\bphi_E(F))$$
In particular, we pick $F=\mo_{e_A}$ to be the structure sheaf of the unit in $A$, then $\bphi_E(F)=E_B$, and if $\Phi(e_A,\alpha)=(b,\beta)$ it follows that 
$$\bphi_E(L_\alpha\otimes F)=L_\beta\otimes T_b^*(\bphi_E(F))$$
$$\Rightarrow \bphi_E(\mo_{e_A})=L_\beta\otimes T_b^*E_B$$
$$\Rightarrow E_B=L_\beta\otimes T_b^*E_B$$
so the image of $\wa$ is in $\{(b,\beta)|(-b,\beta)\in B_{\delta(E_B)}\}$, which is isomorphic to $B_{\delta(E_B)}$. So $\wa$ can be considered as a subabelian variety of $B_{\delta(E_B)}$. And this two have the same dimension, so $\wa\cong B_{\delta(E_B)}$
\end{proof}

We consider the case if $A$ admits a principal polarization. The following theorem is first proved by Orlov \cite{orlov2002derived}, we restate the proof here.

\begin{theorem}
Let $A$ be a principally polarized abelian variety with Picard number 1, then any derived equivalent abelian variety with $A$ is isomorphic to $A$.
\end{theorem}
\begin{proof}
Let $B$ be an abelian variety that is derived equivalent to $A$. Then by theorem \ref{essential}, $\wb\cong A_{\mu}$ for some slope $\mu$. But now we assume $A$ has Picard number $1$, so $\bn(A)\cong\bz$,we choose a generator $L$. Then since $A$ is principally polarized, so we have $deg(L)=1$. Let $\mu=\frac{L^{\otimes n}}{l}$. To consider $A_\mu$, we may assume $(n,l)=1$. Let's consider the map $f:A\to A_\mu$ defined by $a\to(la,h_{L^{\otimes n}}(a))$. So $Ker(f)=A_l\cap K(L^{\otimes n})$. Since we are in the case of char 0, so all these group schemes are reduced. Then we have $A_l$ has $l^{2g}$ elements where $g=\dim A$, and $K(L^{\otimes})$ has $deg(L^{\otimes n})^2=n^{2g}$, but $(n,l)=1$, so $Ker(f)=A_l\cap K(L^{\otimes n})=\{e\}$. So $A\cong A_\mu\cong \wb$. Then $B\cong\wa\cong A$ since $A$ is principally polarized. 
\end{proof}

\section{Semi-homogeneous Vector Bundles as Fourier-Mukai Kernels}

By theorem \ref{Iso-abelian} for any two derived equivalent abelian varieties we can choose a simple semi-homogeneous vector bundle $E$ on $A\times B$ such that $\bphi_E$ is a \fm equivalence. So $\bphi_E(\mo_a)$ is a semi-homogeneous vector bundle. In this section we will consider a \fm equivalence with kernel to be a vector bundle. We will get the isomorphism $A\times \wa\to B\times\wb$ directly.

We need the following lemma in this section.

\begin{lem}
Let $E$ be a simple semi-homogeneous vector bundle on $A$ with slope $\mu$. Then let $\pi_1:A_\mu\to A$ and $\pi_2:A_\mu\to\wa$ be two natural morphisms induced by the projections of $A\times\wa$ to each factor. The $Ker(\pi_1)=\Sigma(E)$ and $deg(\pi_1)=rank(E)^2$.
\end{lem}
\begin{proof}
See Mukai \cite{mukai1978semi}.
\end{proof}

Let $E$ be a simple semi-homogeneous vector bundle on $A\times B$.
\begin{prop}\label{1}
The \fm transform $\bphi_E$ is a \fm equivalence if and only if $E|_{a\times B}$ is simple for some (hence any) point $a\in A$ and for different $a,a'\in A$, we have $E|_{a\times B}$ is not isomorphic to $E|_{a'\times B}$.
\end{prop}
\begin{proof}
By theorem \ref{equiv} and proposition  \ref{im}, we just need to show $E|_{a\times B}$ and $E|_{a'\times B}$ have the same slope.

Let $x=a'-a$, then $T_{(x,e_B)}^*E=E\otimes P_{(\alpha,\beta)}$. So we have $$E|_{a'\times B}\cong E|_{a\times B}\otimes P_\beta$$
so $$\bd (E|_{a'\times B})=\bd(E|_{a\times B})\otimes P_\beta^{\otimes rank(E)}$$
represent the same element in $\bn(B)$, so these two simple semi-homogeneous vectors bundles have the same slope.
\end{proof}

This is a condition that actually hard to control. Next we do some calculation.

Let $E$ on $A\times B$ defines an equivalence. Denote $E_B=E|_{e_A\times B}$ and $E_A=E|_{A\times e_B}$. Denote $M=\bd(E)$ and $M_A,M_B$ similarly. Let $q_1,q_2$ be two projections from $A\times B$ to two factors, and $l=rank(E)$.

Consider $L=M\otimes q_1^*M_A^{-1}\otimes q_2^*M_B^{-1}$. Then $L$ is a line bundle on $A\times B$ satisfies: (1) $L|_{A\times e_B}\cong\mo_A$. (2) For $a\in A$,  $L|_{a\times B}\in\bP^0(B)$. So by the universal property of the Poincare bundle $\mP$ on $\wb\times B$, there exists a morphism $\pi:A\to\wb$ such that $L\cong(\pi\times 1)^*\mP$. 

\begin{lem}
The morphism $\pi$ is an isogeny.
\end{lem}
\begin{proof}
We need to prove $Ker(\pi)$ is finite. Let $\pi(a)=e_{\wb}$, then $L|_{a\times B}\cong\mo_B$, that means if $E|_{a\times B}=E|_{e_A\times B}\otimes P_\beta$ for some $P_\beta\in\bP^0(B)$ (see the proof of the above proposition), then $P_\beta^{\otimes l}=\mo_B$, so $\beta\in\wb_l$. Now different $a\in Ker(\pi)$ have different $\beta\in\wb$ by proposition \ref{1}, so $|Ker(\pi)|\leq|\wb_l|=l^{2g}$, so $\pi$ is an isogeny.
\end{proof}

So we have $M\cong(\pi\times1)^*\mP\otimes q_1^*M_A\otimes q_2^*M_B$. Denote $\mu=\frac{M}{l}$ be the slope of $E$, let's describe $(A\times B)_\mu$. In the following calculation we will use $0$ standing for the unit of any abelian variety if no confusion.

We need to calculate $$T_{(a,b)}^*M\otimes M^{-1}=(T_{(a,b)}^*(\pi\times1)^*\mP\otimes(\pi\times1)^*\mP^{-1})\otimes q_1^*(T_a^*M_A\otimes M_A^{-1})\otimes q_2^*(T_b^*M_B\otimes M_B^{-1})$$
$$\cong q_1^*P_{\hat{\pi}(b)}\otimes q_2^*L_{\pi(a)}\otimes q_1^*(T_a^*M_A\otimes M_A^{-1})\otimes q_2^*(T_b^*M_B\otimes M_B^{-1})$$

here we use $P_\alpha$ (resp. $L_\beta$) to represent the line bundles on $A$ (resp. $B$) corresponding the point $\alpha\in\wa$ (resp. $\beta\in\wb$). $\hat{\pi}$ is the dual of the isogeny $\pi$. So $$(A\times B)_\mu=\{(la,lb,h_{M_A}(a)+\hat{\pi}(b),h_{M_B}(b)+\pi(a))|(a,b)\in A\times B\}$$ (for the definition of $h_L$ see section 3). If we denote $f:(A\times B)_\mu\to A\times B$ induced by the canonical projection to $A\times B$, then $$|Ker(f)|=l^2$$ by the lemma. 

Let $b=0\in B$, then by the theorem \ref{semi}, $T_{(la,0)}^*E\cong E\otimes q_1^*P_{h_{M_A}(a)}\otimes q_2^*L_{\pi(a)}$. Restrict to $0\times B$ we have $E|_{la\times B}\cong E_B\otimes L_{\pi(a)}$. Then since $\bphi_E$ is an equivalence by proposition \ref{1}, we have $E_B\otimes L_{\pi(a)}\cong E_B$ if and only $la=0$, i.e. $a\in A_l$. So we have $$|Ker(\pi)|=\frac{l^{2g}}{|\Sigma(E_B)|}=\frac{l^{2g}}{l^2}=l^{2g-2}$$
so we get the degree of $\pi$.

Next we need to consider $Ker(f)$. Let $g:A\times B\to(A\times B)_\mu$ be the morphism. Recall $$Ker(f)=g(A_l\times B_l)$$
so $$Ker(f)=\{(h_{M_A}(a)+\hat{\pi}(b),h_{M_B}(b)+\pi(a))\in\wa\times\wb~ with ~(a,b)\in A_l\times B_l\}$$
but we know $deg(\pi)=l^{2g-2}$, so $|\{\pi(a)~with~a\in A_l\}|=l^2$. This gives us
$$l^2=|Ker(f)|$$
$$=|\{(h_{M_A}(a)+\hat{\pi}(b),h_{M_B}(b)+\pi(a))\in\wa\times\wb~ with ~(a,b)\in A_l\times B_l\}|$$
$$\geq|\{(h_{M_A}(a),\pi(a))\in\wa\times\wb~ with ~a\in  A_l\}|$$
$$\geq|\{\pi(a)~with~a\in A_l\}|$$
$$=l^2$$
so all $"\geq"$ should be $"="$. In particular, from the first $"\geq"$, we get
$$\{(h_{M_A}(a),\pi(a))\in\wa\times\wb~ with ~a\in  A_l\}=\{\{(\hat{\pi}(b),h_{M_B}(b))\in\wa\times\wb~ with ~b\in  B_l\}\}$$
By the second inequality, we get
$$\{(h_{M_A}(a),\pi(a))\in\wa\times\wb~ with ~a\in  A_l\}=\{\pi(a)~with~a\in A_l\}$$
This tells us $$Ker(\pi)\subset Ker(h_{M_A})\cap A_l$$
Similar
$$Ker(\hat{\pi})\subset Ker(h_{M_B})\cap B_l$$
By using these two equations, let's prove
\begin{theorem}\label{iso}
Notations are as above. If $\bphi_E$ is an equivalence, then $(A\times B)_\mu$ gives an isomorphism between $A\times \wa$ and $B\times \wb$
\end{theorem}
\begin{proof}
We need to prove $p:(A\times B)_\mu\to B\times \wb$ is an isomorphism. We just need to show $Ker(p)=0$. Suppose $(la,lb,h_{M_A}(a)+\hat{\pi}(b),h_{M_B}(b)+\pi(a))\in Ker(p)$, then by definition we have $$(lb,h_{M_B}(b)+\pi(a))=(0,0)\in B\times\wb$$
So $b\in B_l$. So by proposition \ref{semi}
$$T_{(la,0)}^*E\cong E\otimes q_1^*P_{h_{M_A}(a)+\hat{\pi}(b)}$$
restrict to $0\times B$ we get
$$E|_{la\times B}\cong E_B$$
since $\bphi_E$ is an equivalence, by proposition \ref{1} we have $la=0$, i.e. $a\in A_l$. Then by the above discussion we can write (since $(a,b)\in A_l\times B_l$)
$$(h_{M_A}(a)+\hat{\pi}(b),h_{M_B}(b)+\pi(a))=(h_{M_A}(a'),\pi(a'))$$
for some $a'\in A_l$. Now $h_{M_B}(b)+\pi(a)=0$ means $\pi(a')=0$. But $Ker(\pi)\subset Ker(h_{M_A})\cap A_l$, so $h_{M_A}(a')=0$, this proves $$Ker(p)=0$$ so we have an isomorphism $$(A\times B)_\mu\cong B\times \wb$$ similar $$(A\times B)_\mu\cong A\times\wa$$ this gives an isomorphism
$$\eta:A\times \wa\to B\times\wb$$
\end{proof}
Under this isomorphism we can see $$\eta(\wa)=B_{\delta(E_B)}$$ Note that $\eta$ is not an isometry. This isomorphism is a little different with the isomorphism defined in theorem \ref{isometry}. By this isomorphism, the map $\pi$ we defined above correspond to the natural ones $\widehat{B_\mu}\to\wb$ (dual of $B\to B_\mu$)

We can see if we have two abelian varieties $A,B$ and $\bphi_P:D^b(A)\to D^b(B)$ is an equivalence with $\bphi_P(\mo_{e_A})=E_B$ is a simple semi-homogeneous vector bundle on $B$ with rank $l$, then $\wa\cong B_{\delta(E_B)}$. Now we consider the converse problem: For which semi-homogeneous vector bundle $E_B$ on $B$, there exists an abelian variety $A$ and a \fm equivalence $\bphi_P$ such that $\bphi_P(e_A)=E_B$? We can see in this case $P$ must be a semi-simple homogeneous vector bundle on $A\times B$. Here if we have an ample line bundle $K$ on $B$, by replacing $E_B$ with $E_B\otimes K^{\otimes v}$ for $v$ large enough we may assume $M_B$ is ample. Since tensor with a line bundle is also an auto-equivalence, so for the existence of the equivalence, we just assume $\bd(E_B)$ is ample. For all the discussions below, $M_B=\bd(E_B)$ will be ample.

We first see if such $A$ exist, then $A\cong\widehat{B_\mu}$. Similar let $M_B=\bd(E_B)$, we have the following sufficient condition
\begin{lem}
With notations as above, and denote $\pi:A\to\wb$ be the dual of the natural morphism $B\to B_\mu$. Then if there exists a line bundle $M_A$ on $A$ such that the simple semi-homogeneous vector bundle of the slope $$\mu=\frac{q_1^*M_A\otimes(\pi\times1)^*\mP\otimes q_2^*M_B}{l}$$
has rank $l$, then there is a simple semi-homogeneous vector bundle $E$ on $A\times B$ such that $\bphi_E$ is an equivalence and $\bphi_E(\mo_{e_A})=E_B$.
\end{lem}
\begin{proof}
Let $E$ be a simple semi-homogeneous vector bundle with slope $\mu$. Then similar to the above discussion
$$(A\times B)_\mu=\{(la,lb,h_{M_A}(a)+\hat{\pi}(b),h_{M_B}(b)+\pi(a))|(a,b)\in A\times B\}$$
For every $a\in A,b\in B$ we have
$$T_{(la,lb)}^*E\cong E\otimes q_1^*P_{h_{M_A}(a)+\hat{\pi}(b)}\otimes q_2^*L_{h_{M_B}(b)+\pi(a)}$$

Since $E$ is of rank $l$, so $\bd(E)=q_1^*M_A\otimes(\pi\times1)^*\mP\otimes q_2^*M_B$ up to $\bP^0(A\times B)$. That means
$$\bd(E|_{e_A\times B})=\bd(E)|_{e_A\times B}$$
$$=q_1^*M_A\otimes(\pi\times1)^*\mP\otimes q_2^*M_B|_{e_A\times B}$$
$$\cong M_B$$
up to $\bP^0(B)$. So $E|_{e_A\times B}$ has the same slope as $E_B$. Since $E_B$ is simple, we may assume $E|_{e_A\times B}\cong E_B$.

So
$$E|_{-la\times B}\cong E_B\otimes L_{\pi(a)}$$
By the definition of $\pi$, if $l(a-a')\neq 0$, then $\pi(a-a')\neq 0$, so by proposition \ref{1}, $\bphi_E$ is an equivalence. We finish the proof.

\end{proof}


Then let's consider in which case such $N$ exists. We have the following:

\begin{lem}
If there exists a line bundle $N$ on $A$ such that $K(N)\cap A_l=Ker(\pi)$, then this $N$ satisfies the condition in the previous theorem. 
\end{lem}
\begin{proof}
Consider the morphism $\phi:A\to A\times\wa$ defined by $a\to(la,h_N(a))$. Then $\phi$ factors through $\pi$. To check the condition in the previous lemma, we just need to show
$$\{(h_{N}(a)+\hat{\pi}(b),h_{M_B}(b)+\pi(a))|(a,b)\in A_l\times B_l\}$$
has cardinal $l^2$, this is equivalent to 
$$\{(h_{N}(a),\pi(a))|a\in A_l\}=\{(\hat{\pi}(b),h_{M_B}(b))|b\in B_l\}$$

Since $K(N)\cap A_l=Ker(\pi)$ so $h_N$ factors through $\pi$, let $h_N=\theta\pi$. Then $$\{(h_N(a),\pi(a))|a\in A_l\}=\{(\theta\pi(a),\pi(a))|a\in A_l\}$$
and $Ker(h_N)\cap A_l=Ker(\pi)$, so we can check easily 
$$\{(h_{N}(a),\pi(a))|a\in A_l\}=\{(\hat{\pi}(b),h_{M_B}(b))|b\in B_l\}$$
so $N$ satisfies the condition.
\end{proof}

\begin{rem}
Put $N=\bd(E_A)$, we can see the previous lemma is also necessary. 
\end{rem}

We collect the facts we prove in this section here .

\begin{theorem}\label{main}
Let $A,B$ be two abelian varieties. 

(a) If there exists a simple semi-homogeneous vector bundle $E$ on $A\times B$ such that the \fm transform $\bphi_E$ is an equivalence, then:

(1) $\wa\cong B_{\delta(E|_{e_A\times B})}$.

(2) Let $\mu$ be the slope of $E$, then $(A\times B)_\mu$ gives an isomorphism.

(3) Let $A\cong\hat{B_{\delta(E|_{e_A\times B})}}\to\wb$ be the canonical morphism, then $$\bd(E)\cong q_1^*\bd(E|_{A\times e_B})\otimes(\pi\times 1)^*\mP\otimes q_2^*\bd(E|_{e_A\times B})$$

(b) Let $E_B$ be a simple semi-homogeneous vector bundle on $B$ with rank $l$ and slope $\delta$. Let $\pi:B\to B_\delta$ be the natural map. The necessary and sufficient condition for the existence of $\bphi_E$ with $\bphi_E(\mo_{e_A})=E_B$ is there exists a line bundle $N$ on $\hat{B_\delta}$ such that $K(N)\cap\hat{B_\delta}_l=Ker(\hat{\pi})$.
\end{theorem}

\begin{ex}
Let $A$ be an abelian variety, and let $E$ be a simple semi-homogeneous vector bundle on $A$ of rank $l$. If $M=\bd(E)$ satisfies $$deg(M)|l$$ then there exists a semi-homogeneous vector bundle $P$ on $A\times A$ with $\bphi_P$ is an equivalence and $\bphi_P(\mo_{e_A})=E$.
\end{ex}

\bibliographystyle{abbrv}
\bibliography{MyCitation}

\end{document}